# Harmonic Weak Siegel Maaß Forms I
Preimages of Non-Holomorphic Saito-Kurokawa Lifts

Martin Westerholt-Raum

**Abstract:** Given a non-holomorphic Saito-Kurokawa lift we construct a preimage under the vector-valued lowering operator. In analogy with the case of harmonic weak elliptic Maaß forms, this preimage allows for a natural decomposition into a meromorphic and a non-holomorphic part. In this way every harmonic weak Siegel Maaß form gives rise to a Siegel mock modular form.

■ **Siegel mock modular forms**
■ **Harish-Chandra and $(\mathfrak{g}, \mathrm{K})$-modules**
■ **Holomorphic parts of Fourier expansions**

MSC Primary    11F37
MSC Secondary  11F46, 11F70, 22E50

$\mathbf{M}$ORE than 10 years ago Zwegers [Zwe02], and Bruinier and Funke [BF04] independently found definitions of harmonic weak Maaß forms. Zwegers focused on mock theta functions whose modular properties had been a long standing miracle since Ramanujan came up with them in 1920 in his death bed letter. Bruinier and Funke were inspired by the Kudla program and they used harmonic weak Maaß forms to describe dualities between theta lifts. Both approaches to harmonic weak Maaß forms together produced a rich and novel research area, leading to the resolution of vital conjectures in combinatorics, topology, and many other fields—see, for example, [Ono09] for an overview and [BO06; BO10a; DIT11; DMZ12] for some of the applications. In this paper, we suggest an extension of the concept of harmonic weak Maaß forms to the case of Siegel modular forms of genus 2, by building up on the ideas of [BF04].

In order to state some defining formulas, we let $\tau = x + iy \in \mathbb{H}$ be a variable in the Poincaré upper half plane and $0 \le k$ an even integer. Analytic aspects of the theory of harmonic weak Maaß forms have been dominated by two facts. First, the weight-$k$ hyperbolic Laplace operator $\Delta_k$ can be written as a composition $-\xi_{2-k} \circ \xi_k$ of two $\xi$-operators defined by $\xi_k(f) = 2iy^k \overline{\partial_{\overline{\tau}} f}$. Second, the $\xi$-operator gives rise to a short exact sequence of the space $^!\mathrm{M}_k$ of weakly holomorphic modular forms, the space $\mathrm{S}_{2-k}$ of cusp forms, and the space $\mathbb{S}_k$ of harmonic weak Maaß forms whose image under $\xi_k$ is a cusp form.

$$0 \longrightarrow {}^!\mathrm{M}_k \longrightarrow \mathbb{S}_k \xrightarrow{\xi_k} \mathrm{S}_{2-k} \longrightarrow 0$$

This sequence puts the operator $\xi_k$ into the center of the theory. One can define harmonic weak Maaß forms as the modular preimages of weakly holomorphic modular forms of weight $2-k$ under $\xi_k$. Note that the $\xi$-operator is essentially the Maaß lowering operator L. Concretely, we have $\xi_k(f) = y^{2-k} \overline{\mathrm{L} f}$.

Harmonic weak Siegel Maaß forms that we study in this paper are vector-valued. Recall that every complex representation of $\mathrm{GL}_2(\mathbb{C})$ can be viewed as a weight for genus 2 Siegel modular forms. The classical case of weight $k$ modular forms corresponds to the representation $\det^k : \mathrm{GL}_2(\mathbb{C}) \to \mathbb{C}, g \mapsto \det(g)^k$. In our case of genus 2 Siegel modular forms, every complex, irreducible representation of $\mathrm{GL}_2(\mathbb{C})$ can be written as a (tensor) product $\det^k \mathrm{sym}^l$ where $\mathrm{sym}^l$ is the $l$-th symmetric power representation of $\mathrm{GL}_2(\mathbb{C})$. It is necessary to introduce vector-valued Siegel modular forms, because the above target space $\mathrm{S}_{2-k}$ of $\xi_k$ will be replaced by

$${}^{\mathrm{SK}}\mathrm{S}^{(2)}\big(\det^{-k/2}\mathrm{sym}^k\big) : \text{space of non-holomorphic Saito-Kurokawa lifts of weight } k \text{ cusp forms.}$$

We revisit the non-holomorphic Saito-Kurokawa lift in Section 1. Note that non-holomorphic Saito-Kurokawa lifts do not occur for scalar weights $\det^k$ for any $k$.

The space of harmonic weak Siegel Maaß forms is defined as the space of real-anlaytic functions with possible meromorphic singularities that transform like Siegel modular forms and that are mapped to non-holomorphic Saito-Kurokawa lifts under the vector-valued lowering operators L. This space contains the space of meromorphic Siegel modular forms $^!\mathrm{M}^{(2)}\big(\mathrm{det}^{2-k/2}\mathrm{sym}^{k-2}\big)$, which are annihilated by the lowering operator. In analogy with the above short exact sequence for harmonic weak elliptic modular forms, we show that every non-holomorphic Saito-Kurokawa lift allows for a modular preimage.

**Theorem I.** *There is an exact sequence of weak Siegel Maaß forms*

$$0 \longrightarrow {}^!\mathrm{M}^{(2)}\big(\mathrm{det}^{2-k/2}\mathrm{sym}^{k-2}\big) \longrightarrow {}^{\mathrm{SK}}\mathbb{S}^{(2)}\big(\mathrm{det}^{2-k/2}\mathrm{sym}^{k-2}\big) \xrightarrow{\mathrm{L}} {}^{\mathrm{SK}}\mathrm{S}^{(2)}\big(\mathrm{det}^{-k/2}\mathrm{sym}^{k}\big) \longrightarrow 0$$

Elements of of $^{\mathrm{SK}}\mathbb{S}^{(2)}\big(\mathrm{det}^{2-k/2}\mathrm{sym}^{k-2}\big)$ will henceforth be called harmonic weak Siegel Maaß forms.

*Remark.* The vector-valued lowering operator defined in [Kle+15], which we will use in this paper takes weight-$\mathrm{det}^k\mathrm{sym}^l$ forms to forms of weight $\mathrm{det}^{k-2}\mathrm{sym}^2\mathrm{sym}^l$, containing the weight $\mathrm{det}^{k-2}\mathrm{sym}^{2+l}$ as a subrepresentation. The above sequence should be read in such a way that the image of harmonic weak Siegel Maaß forms under L is a priori modular of weight $\mathrm{det}^{-k/2}\mathrm{sym}^2\mathrm{sym}^{k-2}$ and that this image vanishes outside of $\mathrm{det}^{-k/2}\mathrm{sym}^k \subset \mathrm{det}^{-k/2}\mathrm{sym}^2\mathrm{sym}^{k-2}$. This in particular, implies that the kernel vanishes under L, which characterizes meromorphic Siegel modular forms.

*Remark.* The above kind of harmonic weak Siegel Maaß forms is by no means the only one that we expect to exist. It is so far unclear how to construct harmonic weak Siegel Maaß forms whose image under the raising operator is a Siegel holomorphic modular form. This space would naturally contain Siegel modular forms with singularities whose analytic behavior is analogous with the one of non-holomorphic Saito-Kurokawa lifts. Even such Siegel modular forms have not yet been explored.

The occurrence of vector-valued Siegel modular forms cannot be avoided, since non-holomorphic Saito-Kurokawa lifts do not allow for a scalar-valued realization. However, harmonic weak Siegel Maaß forms can be forced into scalar weights by applying the vector valued raising operator sufficiently often. Denote the space of weight $k/2$ almost meromorphic Siegel modular forms of depth $k/2 - 1$ by $^!\mathrm{M}^{(2)}\big((\mathrm{det}^{k/2})^{[k/2-1]}\big)$; A definition and study of almost meromorphic Siegel modular forms can be found in [Kle+15]. We define almost harmonic weak Siegel Maaß forms of weight $\mathrm{det}^{k/2}$ as real-analytic modular forms with singularities whose image under $\mathrm{L}^{k/2}$ is a non-holomorphic Saito-Kurokawa lift. We denote the space of such forms by $^{\mathrm{SK}}\mathbb{S}^{(2)}\big((\mathrm{det}^{k/2})^{[k/2-1]}\big)$.

**Corollary II.** *There is an exact sequence of weak Siegel Maaß forms*

$$0 \longrightarrow {}^!\mathrm{M}^{(2)}\big((\mathrm{det}^{k/2})^{[k/2-1]}\big) \longrightarrow {}^{\mathrm{SK}}\mathbb{S}^{(2)}\big((\mathrm{det}^{k/2})^{[k/2-1]}\big) \xrightarrow{\mathrm{L}^{k/2}} {}^{\mathrm{SK}}\mathrm{S}^{(2)}\big(\mathrm{det}^{-k/2}\mathrm{sym}^k\big) \longrightarrow 0$$

*Remark.* The power of the lowering operator in the previous statement is a tensor power whose target weight is $\mathrm{det}^{-k/2}\mathrm{sym}^k \subset \mathrm{det}^{-k/2}(\mathrm{sym}^2)^{k/2}$. The exact sequence has to be interpreted as the one in Theorem I.

One source of interest in harmonic weak Maaß forms is their Fourier series expansion and mock modular forms associated with them. Writing $e(x) = \exp(2\pi i x)$ for $x \in \mathbb{C}$, any harmonic weak Maaß form $f$ can be naturally decomposed as

$$f(\tau) = \sum_{n \gg -\infty} c^+(f;n)e(n\tau) + \sum_{n \ll \infty} c^-(f;n)\,{}^-a_k^{(1)}(n,y)e(n\tau),$$

where, if $k \neq 1$, then $^-a_k^{(1)}$ is essentially the $W$-Whittaker function if $n < 0$, the $M$-Whittaker function if $n > 0$, and a power of $y = \Im\tau$ if $n = 0$. The first sum is called the holomorphic part of $f$ and the second one is called its non-holomorphic part. The holomorphic part of $f$ is a mock modular form and conversely by definition every mock modular form has a modular completion of the above form. Note that harmonic weak Maaß forms with $c^-(f;n) = 0$ for $n \geq 0$ are best understood. Their image under $\xi_k$ is a cusp form that essentially has Fourier coefficients $c^-(f;n)$. Their Fourier coefficients $c^+(f;n)$ for $n \leq 0$ are not completely understood, but they were successfully related to derivatives of twisted $L$-values in [BO10b].

To formulate the analogue of the above decomposition, we set $e(x) = \exp(2\pi i \operatorname{trace}(x))$ for any complex square matrix $x$. We put the superscript $\pm$ to the left of the analytic part $a(t, y)$ of the Fourier coefficient, reserving the right superscript for annotation with regards the genus (which in our case is 2).

**Theorem III.** *Every $f \in {}^{\mathrm{SK}}\mathbb{S}^{(2)}\bigl(\det^{2-k/2}\operatorname{sym}^{k-2}\bigr)$ has a Fourier series expansion of the form*

$$f(\tau) = \sum_t c^+(f;t)\, e(t\tau) + \sum_{t\ \text{indefinite}} c^-(f;t)\,{}^-a_k^{(2)}(t,y)e(tx),$$

*where*

$$\sum_{t\ \text{indefinite}} c^-(f;t)\,{}^{\mathrm{SK}}a_k^{(2)}(t,y)e(tx) \in {}^{\mathrm{SK}}\mathrm{S}^{(2)}\bigl(\det^{-\frac{k}{2}}\operatorname{sym}^{k-2}\bigr)$$

*for real-analytic functions ${}^-a_k^{(2)}$ and ${}^{\mathrm{SK}}a_k^{(2)}$ defined in Section 4.*

Given the analytic characterization of harmonic weak Siegel Maaß forms, it seems justified to call the first summand in Theorem III a Siegel mock modular form. We given a precise definition in Section 4.3. We plan to investigate the Fourier coefficients of Siegel mock modular forms in a sequel to this paper.

We discuss the proof of Theorem I. The proof generalizes ideas in [BF04], but it requires additional input from the theory of $(\mathfrak{g}, K)$-modules. The following description of $(\mathfrak{g}, K)$-modules associated with some harmonic weak Maaß forms was independently found by Ralf Schmidt in an unpublished note. Further, the author was notified that Schulze-Pillot wrote about $(\mathfrak{g}, K)$-modules associated with harmonic weak Maaß forms already in [SP11]. Recall that in the case of nonpositive, even weight $k$, there are two differential operators that map harmonic weak Maaß forms to weakly holomorphic ones. One of them, the $\xi_k$-operators has already been mentioned, and another one arises from Bol's identity.

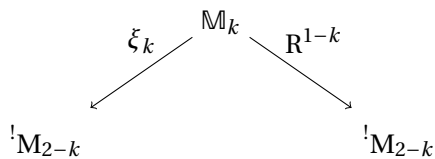

First note that from the perspective of $(\mathfrak{g}, K)$-modules it is not very natural to include complex conjugation into the above diagram. Rather, one sacrifices holomorphicity and employs the lowering operators L instead of $\xi_k$. This yields a map

$$\mathbb{M}_k \longrightarrow y^{2-k}\overline{{}^!\mathrm{M}_{2-k}},$$

where the right hand side consists of functions transforming like modular forms of weight $k-2$. In other words, the $(\mathfrak{g}, K)$-module associated with a non-holomorphic $f \in \mathbb{M}_k$ contains the direct sum of holomorphic and antiholomorphic (limits of) discrete series. The quotient is a finite dimensional and irreducible $(\mathfrak{g}, K)$-module. It corresponds to the images $\mathrm{R}^j f$ with $0 \leq j \leq -k$ of $f \in \mathbb{M}_k$ under powers of the raising operator.

One can phrase the above observation as follows: The $(\mathfrak{g}, K)$-module associated with a proper harmonic weak Maaß form $f$ is isomorphic to a suitable degenerate principal series. This observation guides our proof of existence of harmonic weak Siegel Maaß forms. The Harish-Chandra module attached to non-holomorphic Saito-Kurokawa lifts can be realized as a Langlands quotient. Lee [Lee96] gave a very explicit description of the corresponding principal series. It turns out that they have length 3. The Saito-Kurokawa lift corresponds to the unique irreducible quotient. The second composition factor is the direct sum of a holomorphic and antiholomorphic discrete series, whose minimal $K$-type corresponds to a vector-valued Siegel modular form. The precise description of the Harish-Chandra module delivered by Lee enables us to extend the argument in [BF04] to Siegel modular forms of genus 2.

Bruinier and Funke employ sheaf cohomology. Since they dealt with complex curves, the Dolbeaux resulution of the structure sheaf is surjective onto $(0,1)$-forms. This can be interpreted in such a way that locally for every $(0,1)$-form there is a preimage in 0-forms. Translated into classical language, every weakly holomorphic modular form is locally the image of some harmonic 0-form. This is no longer true in the case of complex 3-folds, which we deal with in this paper. Vanishing under the Dolbeaux derivative $\overline{\partial}$ is neccessary for a $(0,1)$-form to locally admit a preimage under $\overline{\partial}$. We reformulate this vanishing condition in terms of $(\mathfrak{g}, K)$-modules, and establish it for non-holomorphic Saito-Kurokawa lifts.

The main step in the proof of Bruinier and Funke is the second step in ours. To guarantee global existence of preimages of harmonic weak Maaß forms they prove vanishing of a corresponding obstruction space, constituted by the first cohomology of a certain holomorphic line bundle. Using Serre duality they reduce it to a classical vanishing result for elliptic modular forms. In our case, the obstruction space is the first cohomology of a vector bundle. We employ a vanishing theorem by Grauert-Riemenschneider [GR70] and Igusa's study [Igu62] of Siegel modular forms to deduce its vanishing.

We start the paper with preliminaries in Section 1. Section 2 contains the theory of Harish-Chandra modules and degenerate principal series that we need in the paper. We proceed to the proof of existence of harmonic weak Siegel modular forms in Section 3. That section contains most of the complex geometry that is important to the paper. The last Section 4 focuses on Siegel mock modular forms, which arise from a natural splitting of the Fourier expansion of harmonic Siegel modular forms.

**Acknowledgment.** The author is grateful to Özlem Imamoğlu and Olav Richter for many helpful discussions about harmonic weak Maaß forms and their generalizations. He thanks Rainer Weissauer for discussions about principal series and the endoscopy. He appreciates very much comments on an early version of this paper by Kathrin Bringman, Jan H. Bruinier, Olav Richter, and Ralf Schmidt.

# 1 Preliminaries

**§1.1 The symplectic group of genus** 2. Throughout, we focus on the group $G = \mathrm{Sp}_2(\mathbb{R})$ with complex Lie algebra $\mathfrak{g} = \mathfrak{sp}_2$, whose definitions are

$$G = \mathrm{Sp}_2(\mathbb{R}) = \left\{ g = \begin{pmatrix} a & b \\ c & d \end{pmatrix} \in \mathrm{Mat}_4(\mathbb{R}) : {}^\mathrm{t}g J_2 g = J_2 \right\} \quad \text{with} \quad J_2 = \begin{pmatrix} 0 & -I_2 \\ I_2 & 0 \end{pmatrix}, \; I_2 = \begin{pmatrix} 1 & 0 \\ 0 & 1 \end{pmatrix}, \text{ and}$$

$$\mathfrak{g}_0 = \mathfrak{sp}_2(\mathbb{R}) = \left\{ \begin{pmatrix} a & b \\ c & -{}^\mathrm{t}a \end{pmatrix} : a \in \mathrm{Mat}_2(\mathbb{R}), b, c \in \mathrm{Mat}_2^\mathrm{t}(\mathbb{R}) \right\}, \quad \mathfrak{g} = \mathfrak{g}_0 \otimes_\mathbb{R} \mathbb{C}.$$

A Cartan involution $\theta$ of $\mathfrak{g}_\mathbb{R}$ is given by the map sending $\begin{pmatrix} a & b \\ c & -{}^\mathrm{t}a \end{pmatrix}$ to $\begin{pmatrix} -{}^\mathrm{t}a & -c \\ -b & a \end{pmatrix}$. The associated Cartan decomposition is $\mathfrak{g} = \mathfrak{k} \oplus \mathfrak{m}$ with

$$\mathfrak{k} = \left\{ \begin{pmatrix} a & b \\ -b & a \end{pmatrix} : a \in \mathrm{Mat}_2(\mathbb{C}), b \in \mathrm{Mat}_2^\mathrm{t}(\mathbb{C}), a = -{}^\mathrm{t}a \right\} \quad \text{and} \quad \mathfrak{m} = \left\{ \begin{pmatrix} a & b \\ b & -a \end{pmatrix} : a, b \in \mathrm{Mat}_2^\mathrm{t}(\mathbb{C}) \right\}.$$

The corresponding compact subgroup $\mathrm{U}_2(\mathbb{R}) \cong K(\mathbb{R}) \subset G(\mathbb{R})$ consists of matrices $\begin{pmatrix} a & b \\ -b & a \end{pmatrix} \in \mathrm{Sp}_2(\mathbb{R})$ with $a - ib \in \mathrm{U}_2(\mathbb{R})$. We write $\mathfrak{g}_0$ and $\mathfrak{k}_0$ for the real Lie algebras of $G(\mathbb{R})$ and $K(\mathbb{R})$.

We now fix a bases of $\mathfrak{k}$ and of $\mathfrak{m}$. Define $\mathfrak{h}_\mathfrak{c} = -i \begin{pmatrix} & -I_2 \\ I_2 & 0 \end{pmatrix}$ and

$$\mathfrak{h}_\mathfrak{k} = \begin{pmatrix} & -i & & \\ & & i & \\ i & & & \\ & & & -i \end{pmatrix}, \quad \mathfrak{e}_\mathfrak{k} = \frac{1}{2}\begin{pmatrix} & i & & 1 \\ -i & & 1 & \\ & -1 & & i \\ -1 & & -i & \end{pmatrix}, \quad \mathfrak{f}_\mathfrak{k} = \frac{1}{2}\begin{pmatrix} & i & & -1 \\ -i & & -1 & \\ & 1 & & i \\ 1 & & -i & \end{pmatrix},$$

$$\mathfrak{h}_\mathfrak{m}^+ = \begin{pmatrix} -i & & -1 & \\ & -i & & -1 \\ -1 & & i & \\ & -1 & & i \end{pmatrix}, \quad \mathfrak{e}_\mathfrak{m}^+ = \begin{pmatrix} & & 1 & -i \\ & & -i & -1 \end{pmatrix}, \quad \mathfrak{e}_\mathfrak{m}^+ = \begin{pmatrix} 1 & & -i & \\ -i & & -1 & \end{pmatrix},$$

$$\mathfrak{h}_\mathfrak{m}^- = \begin{pmatrix} i & & -1 & \\ & i & & -1 \\ -1 & & -i & \\ & -1 & & -i \end{pmatrix}, \quad \mathfrak{e}_\mathfrak{m}^- = \begin{pmatrix} 1 & & i & \\ i & & -1 & \end{pmatrix}, \quad \mathfrak{e}_\mathfrak{m}^- = \begin{pmatrix} & & 1 & i \\ & & i & -1 \end{pmatrix}.$$

Names are chosen in such a way that the action of $\mathfrak{k}$ is the customary one: $[\mathfrak{e}_\mathfrak{k}, \mathfrak{e}_\mathfrak{m}^+] = 0$, $[\mathfrak{e}_\mathfrak{k}, \mathfrak{h}_\mathfrak{m}^+] = -2\mathfrak{e}_\mathfrak{m}^+$, etc. The element $\mathfrak{h}_\mathfrak{c}$ acts on $\mathfrak{m}^\pm$ by multiplication with $\pm 2$. Compact and non-compact roots of $\mathfrak{g}$ are

$$\pm 2\mathfrak{h}_\mathfrak{k}^\vee; \quad \pm 2\mathfrak{h}_\mathfrak{c}^\vee, \pm 2(\mathfrak{h}_\mathfrak{c}^\vee + \mathfrak{h}_\mathfrak{k}^\vee), \pm 2(\mathfrak{h}_\mathfrak{c}^\vee - \mathfrak{h}_\mathfrak{k}^\vee).$$

Note that $i\mathfrak{h}_\mathfrak{c}$ and $i\mathfrak{h}_\mathfrak{k}$ span a Cartan subalgebra of both $\mathfrak{k}_0$ and $\mathfrak{g}_0$. Thus we can, as is customary, identify weights of $K$ and $G$.

**§1.2 The upper half space.** The genus 2 Siegel upper half space $\mathbb{H}^{(2)} = \{\tau = x + iy \in \mathrm{Mat}_2^\mathfrak{t}(\mathbb{C}) : y > 0\}$ consists of symmetric complex matrices $\tau$ with positive definite imaginary part. It carries an action of $\mathrm{Sp}_2(\mathbb{R})$ by means of

$$\mathrm{Sp}_2(\mathbb{R}) \circlearrowright \mathbb{H}^{(2)} : \begin{pmatrix} a & b \\ c & d \end{pmatrix} \tau = (a\tau + b)(c\tau + d).$$

The entries of $\tau$, $x$, and $y$ are given by

$$\tau = \begin{pmatrix} \tau_1 & z \\ z & \tau_2 \end{pmatrix}, \quad y = \begin{pmatrix} y_1 & v \\ v & y_2 \end{pmatrix}, \text{ and } \quad x = \begin{pmatrix} x_1 & u \\ u & x_2 \end{pmatrix}.$$

We associate to every $\tau \in \mathbb{H}^{(2)}$ the element

$$g_\tau = \begin{pmatrix} \sqrt{y} & x\sqrt{y}^{-1} \\ & \sqrt{y}^{-1} \end{pmatrix} \in \mathrm{Sp}_2(\mathbb{R}),$$

where $\sqrt{y}$ is the unique positive definite symmetric square root of $y$.

**§1.3 Weights, types, and slash actions.** We call a finite dimensional, complex representation $\sigma$ of $\mathrm{GL}_2(\mathbb{C})$ a weight for genus 2 Siegel modular forms. In genus 2, any weight is isomorphic to a direct sum of representations $\det^k \mathrm{sym}^l$ for some $k \in \mathbb{Z}$ and some $0 \le l \in \mathbb{Z}$. The representation space of $\sigma$ will be denoted by $V(\sigma)$. A finite dimensional, complex representation of $\Gamma^{(2)} = \mathrm{Sp}_2(\mathbb{Z})$ is called a type. We focus throughout and without further mentioning it on types whose kernel has finite index in $\Gamma^{(2)}$. The representation space of $\rho$, in analogy with weights, is denoted by $V(\rho)$. By slight abuse of notation we denote the trivial weight and trivial type by the same letter $\mathbb{1}$.

A $\mathrm{GL}_2(\mathbb{C})$-valued cocycle is given by $j(g, \tau) = c\tau + d$ for $g \in \mathrm{Sp}_2(\mathbb{R})$ and $\tau \in \mathbb{H}^{(2)}$. Fixing a weight $\sigma$ and a type $\rho$, we define the slash action of weight $\sigma$ and type $\rho$ on functions $f : \mathbb{H}^{(2)} \to V(\sigma) \otimes V(\rho)$ by

$$f|_{\sigma, \rho}\, g = \sigma(c\tau + d)^{-1} \otimes \rho(\gamma)^{-1} (f \circ g).$$

If the type is trivial, we suppress it in the notation, writing $|_\sigma$ instead of $|_{\sigma, \rho}$.

**§1.4 Lowering and raising operators.** Order 1 lowering and raising operators for Siegel modular forms are necessarily vector-valued. Existence and basic properties follow from general work by Helgason [Hel59; Hel62]. Concrete expressions in classical terms were computed in [Kle+15]. As usual, we define the matrix differentials

$$\partial_\tau = \begin{pmatrix} \partial_{\tau_1} & \tfrac{1}{2}\partial_z \\ \tfrac{1}{2}\partial_z & \partial_{\tau_2} \end{pmatrix} \quad\text{and}\quad \partial_{\overline\tau} = \begin{pmatrix} \partial_{\overline{\tau_1}} & \tfrac{1}{2}\partial_{\overline z} \\ \tfrac{1}{2}\partial_{\overline z} & \partial_{\overline{\tau_2}} \end{pmatrix}.$$

The lowering operator can be written independently of the weight. Given a differentiable function $f : \mathbb{H}^{(2)} \to V(\sigma) \otimes V(\rho)$, we set

$$\mathrm{L} f = \mathrm{L}_\sigma f = y^{\mathrm{t}}\bigl( y \partial_{\overline\tau} \bigr) \otimes f.$$

The raising operator depends on the weight, and we give a formula specific to weights $\det^k\mathrm{sym}^l$. We define a symmetrization map that was devised in [Kle+15]. Observe that elements of $V(\mathrm{sym}^l)$ can be written as a (tensor) product of $l$ elements in the standard representation std of $\mathrm{GL}_2(\mathbb{C})$, which satisfies $V(\mathrm{std}) \cong \mathbb{C}^2$. For $0 \le l \in \mathbb{Z}$, we define

$$\mathrm{t}_{\det^k\mathrm{sym}^l} : \mathrm{sym}^2 \otimes \det^k\mathrm{sym}^l \longrightarrow \mathrm{sym}^2 \otimes \det^k\mathrm{sym}^l,$$

$$v_1 v_2 \otimes \prod_{j=1}^{l} v_{2+j} \longmapsto \frac{1}{2l} \sum_{j'=1}^{l} \Bigl( v_{2+l} v_2 \otimes v_1 \prod_{\substack{j=1\\ j\ne j'}}^{l} v_{2+j} + v_1 v_{2+l} \otimes v_2 \prod_{\substack{j=1\\ j\ne j'}}^{l} v_{2+j} \Bigr).$$

The raising operator on a differentiable function $f : \mathbb{H}^{(2)} \to V(\det^k\mathrm{sym}^l) \otimes V(\rho)$ is given by

$$\mathrm{R} f = \mathrm{R}_{\det^k\mathrm{sym}^l} f = \partial_\tau \otimes f - \frac{ik}{2} y^{-1} \otimes f - \frac{il}{2} \mathrm{t}_{\det^k\mathrm{sym}^l}\bigl( y^{-1} \otimes f \bigr).$$

We will suppress the weight from notation, if it becomes clear form the context.

The covariance properties of L and R define them uniquely up to scalar multiples by Proposition 6.9 of [Kle+15]. For every $g \in \mathrm{Sp}_2(\mathbb{R})$, we have

$$(\mathrm{L} f)\big|_{\mathrm{sym}^{\vee 2}\sigma}\, g = \mathrm{L}\bigl( f\big|_\sigma g \bigr) \quad\text{and}\quad (\mathrm{R} f)\big|_{\mathrm{sym}^2\sigma}\, g = \mathrm{R}\bigl( f\big|_\sigma g \bigr), \tag{1.1}$$

where $\mathrm{sym}^{\vee 2}$ denotes the dual representation of $\mathrm{sym}^2$. To simplify stating the next definition, we set $\mathrm{L}\sigma = \mathrm{sym}^{\vee 2}\sigma = \det^{-2}\mathrm{sym}^2\sigma$ and $\mathrm{R}\sigma = \mathrm{sym}^2\sigma$. By viewing $\mathrm{GL}_2(\mathbb{C})$-representations as $\mathrm{U}_2(\mathbb{R})$-representations, we define the following orthogonal projections for $\sigma = \det^k\mathrm{sym}^l$ if $l \ge 2$.

$$\pi_{\mathrm{L},+} : \mathrm{L}\sigma \to \det^{k-2}\mathrm{sym}^{l+2}, \quad \pi_{\mathrm{L},0} : \mathrm{L}\sigma \to \det^{k-1}\mathrm{sym}^l, \quad \pi_{\mathrm{L},-} : \mathrm{L}\sigma \to \det^k\mathrm{sym}^{l-2};$$

$$\pi_{\mathrm{R},+} : \mathrm{R}\sigma \to \det^k\mathrm{sym}^{l+2}, \quad \pi_{\mathrm{R},0} : \mathrm{R}\sigma \to \det^{k+1}\mathrm{sym}^l, \quad \pi_{\mathrm{R},-} : \mathrm{R}\sigma \to \det^{k+2}\mathrm{sym}^{l-2}.$$

It is clear that $\pi_{\mathrm{L},+} \circ \mathrm{L}$ and all analogue differential operators satisfy covariance properties derived from (1.1).

**§1.5 The non-holomorphic Saito-Kurokawa lift.** Existance of the non-holomorphic Saito-Kurokawa lift was first indicated in work of Howe and Piatetski-Shapiro [HPS83]. Piatetski-Shapiro [PS83] used it to determine cuspical automorphic representations for $\mathrm{PGSp}_2$ whose $L$-function has a pole. A special case of the Saito-Kurokawa lift that is of particular interest from a geometrical perspective, was investigated by Bruinier; see Example 5.13 of [Bru02]. Since the (non-holomorphic) Saito-Kurokawa lift can be expressed as a theta lift (or Weil lifting as Piatetski-Shapiro calls it), it appears in Li's study of theta lifts [Li90]. In the present paper, we will focus on the formulation by Schmidt [Sch05], who established functoriality, and Miyazaki [Miy04], who computed Fourier expansions of non-holomorphic Saito-Kurokawa lifts. Given a real-analytic Siegel modular eigenform $f$ of weight $\det^{-k/2}\mathrm{sym}^k$ with $k \in 2\mathbb{Z}$ for some finite index subgroup $\Gamma' \subseteq \Gamma^{(2)}$, we say that $f$ is a non-holomorphic Saito-Kurokawa lift if its standard $L$-function has a pole. The linear span of such Saito-Kurokawa lifts is denoted by $^{\mathrm{SK}}\mathrm{S}^{(2)}\bigl(\det^{-k/2}\mathrm{sym}^k, \Gamma'\bigr)$.

**Theorem 1.1 (Piatetski-Shapiro).** *Given an elliptic modular form $f$ of level $1$ and weight $k$ with $4 \mid k$, then there is a cuspidal, real-analytic Siegel modular form $g$ of level $1$ and weight $\det^{-k/2} \mathrm{sym}^k$ whose standard L-function has a pole and which contains the L-function of $f$ as a factor. In particular, for $k \geq 12$ divisible by $4$, the space $^{\mathrm{SK}}\mathrm{S}^{(2)}(\det^{-k/2}\mathrm{sym}^k, \Gamma^{(2)})$ is not empty.*

Since we work with vector-valued Siegel modular forms for $\Gamma^{(2)}$, we now explain how to pass from Saito-Kurokawa lifts in the classical sense to vector-valued ones. The reader who prefers to work in a slightly more classical setting can skip this and assume $4 \mid k$ and $\rho = \mathbb{1}$ in the rest of the paper. None of the steps of our construction involves the type in a crucial way.

For a finite index subgroup $\Gamma' \subseteq \Gamma^{(2)}$, we denoted the induction of the trivial representation from $\Gamma'$ to $\Gamma^{(2)}$ by $\rho_{\Gamma'}$. Recall the induction map from modular forms for a subgroup $\Gamma' \subseteq \Gamma^{(2)}$ to modular forms of type $\rho_{\Gamma'}$, which, for example, was spelled out in [WR14]. The image of $^{\mathrm{SK}}\mathrm{S}^{(2)}(\det^{-k/2}\mathrm{sym}^k, \Gamma')$ under the induction map is denoted by $^{\mathrm{SK}}\mathrm{S}^{(2)}(\det^{-k/2}\mathrm{sym}^k \otimes \rho_{\Gamma'})$. The decomposition $\bigoplus_j \rho_j$ of $\rho_{\Gamma'}$ into irreducible components yields a natural decomposition

$$^{\mathrm{SK}}\mathrm{S}^{(2)}(\det^{-k/2}\mathrm{sym}^k \otimes \rho_{\Gamma'}) = \bigoplus_j {}^{\mathrm{SK}}\mathrm{S}^{(2)}(\det^{-k/2}\mathrm{sym}^k \otimes \rho_j),$$

which we employ to define the spaces on the right hand side.

Consider Theorem 5.3 of [Miy04]. It tells us that there is a polynomial $p_k(t, \sqrt{y})$ that takes values in $V(\det^{-k/2}\mathrm{sym}^k)$, and that depends on the entries of indefinite $t \in \mathrm{Mat}_2(\mathbb{R})$ and the positive definite square root $\sqrt{y}$ of $y$ such that any level $1$ Saito-Kurokawa lift has Fourier expansion

$$\sum_{t \text{ indefinite}} c(f; t)\, {}^{\mathrm{SK}}\tilde{a}^{(2)}_{\det^{-k/2}\mathrm{sym}^k}(t; y) e(tx) \tag{1.2}$$

with coefficients $c(f; t) \in \mathbb{C}$ and

$$^{\mathrm{SK}}\tilde{a}^{(2)}_{\det^{-k/2}\mathrm{sym}^k}(t, y) = \det(y)^{\frac{k+2}{4}} p_k(t, \sqrt{y}) \frac{K_{k/2}\big(\mathrm{trace}(ty)^2 - \det(t)\det(y)\big)}{\big(\mathrm{trace}(ty)^2 - \det(t)\det(y)\big)^{k/4}}, \tag{1.3}$$

where $K_{k/2}$ is the $K$-Bessel function. The precise shape of the polynomials $p_k$ is irrelevant to use, since later in Section 4 we replace it by a more suitable choice. Note that the coefficients $c(f; t)$ grow polynomially in the entries of $t$.

From the construction of Saito-Kurokawa lifts as theta lifts it becomes clear that for Saito-Kurokawa lifts of nontrivial type $f \in {}^{\mathrm{SK}}\mathrm{S}^{(2)}(\det^{-k/2}\mathrm{sym}^k \otimes \rho)$ we have an analogous Fourier expansion with $c(f; t) \in V(\rho)$.

## 2 Harish-Chandra modules

This section contains most of the theory of $(\mathfrak{g}, \mathrm{K})$-modules that we apply in the present paper. We outline some of the basic concepts in Subsection 2.1, but generally the reader is referred to [Kna01; Wal88] for details and precise definitions. In Subsection 2.2 we describe how to pass from functions on the upper half space to $(\mathfrak{g}, \mathrm{K})$-modules and how the behavior under differential operators is connected to $K$-types.

Subsection 2.3 contains a description of some Vogan-Zuckerman modules. In Subsection 2.4, we identify the non-holomorphic Saito-Kurokawa lift at the infinite place as a suitable module $A_{\mathfrak{q}}(\lambda)$. In Subsection 2.5, we consider the composition series of some degenerate generalized principal series. In particular, we recognize the Saito-Kurokawa $(\mathfrak{g}, K)$-module as a submodule in a specific principle series.

**§2.1 Preliminaries on Harish-Chandra modules.** Given a Lie group $G$ with complex Lie algebra $\mathfrak{g}$ and a maximal compact subgroup $K \subseteq G$, a $(\mathfrak{g}, \mathrm{K})$-module is a $\mathfrak{g}$-module which carries a compatible action of $K$. The definition is spelled out in Section 3.3.1 of [Wal88]. A $(\mathfrak{g}, \mathrm{K})$-module is called admissible if it is a unitary $K$-module and the $K$-isotypical components are finite dimensional. An admissible $(\mathfrak{g}, \mathrm{K})$-module is called

Harish-Chandra module. It is a theorem by Harish-Chandra (cf. Theorem 3.4.10 in [Wal88]) that every irreducible unitary representation $\pi$ of $G$ gives rise to an admissible $(\mathfrak{g},K)$-module by taking smooth vectors. It is called the Harish-Chandra module associated with $\pi$.

The $K$-types that occur in a $(\mathfrak{g},K)$-module and transitions between them are the only feature that we will make use of. It is not necessary but convenient to draw them. As throughout the whole paper, we focus on the case $G = \mathrm{Sp}_2(\mathbb{R})$. Then irreducible representations of $K$ can be parametrized by dominant weights of $K$, which can be expressed in terms of roots of $\mathfrak{g}$. We denote the weight $(a+b)\mathfrak{h}_{\mathfrak{c}}^{\vee} + (a-b)\mathfrak{h}_{\mathfrak{k}}^{\vee}$ by the pair $(a,b)$. As a $\mathrm{GL}_2(\mathbb{C})$ representation it corresponds to $\det^b \mathrm{sym}^{a-b}$. Possible $K$-types range in $\Lambda = \{(a,b) \in \mathbb{Z}^2 : a \geq b\}$. We depict a $(\mathfrak{g},K)$-module by an arrangement of circles that correspond to all possible $(a,b)$. Filled circles indicate $K$-types that occur with multiplicity at least one, and the remaining ones indicate those that do not occur. The two diagrams below are typical such diagrams of $K$-types. We read the diagram as if it would continue infinitely to the right and to the bottom without indicating this. In the left diagram we have encircled the $K$-isotrivial component, which corresponds to the pair $(0,0)$. In general we will assume that our diagram is symmetric and therefore refrain from marking $(0,0)$.

The action of $\mathfrak{m} \subset \mathfrak{g}$ gives rise to transitions of $K$-types, whose possible targets can be computed by the Clebsch-Gordan rules. For $a-b \geq 2$ the image of a $K$-type $(a,b)$ under $\mathfrak{m}^{\pm}$ is a direct sum of $K$-types $(a\pm 2, b)$, $(a\pm 1, b\pm 1)$, and $(a, b\pm 2)$. The image of $(a+1, a)$ under $\mathfrak{m}^+$ and $\mathfrak{m}^-$ is a direct sum of $(a+3, a)$ and $(a+2, a+1)$, and $(a+1, a-2)$ and $(a, a-1)$, respectively. The image of $(a,a)$ under $\mathfrak{m}^+$ and $\mathfrak{m}^-$ is $(a+2, b)$ and $(a, a-2)$.

For a given $(\mathfrak{g},K)$-module some of those transition maps might be nonzero. This is obviously the case if they hit a unfilled circle, and no further indication in our diagrams is needed. A second phenomenon, which we call (transition) walls, can occur. A wall is a linear relation $w_a a + w_b b = w_c$ for $(a,b)$ and a sign $w_{\sigma} \in \{\pm 1\}$ such that every transition from $(a,b)$ to $(a',b')$ with $w_a a + w_b b = w_c$ is zero, if $w_{\sigma} \cdot (w_a a' + w_b b' - w_c) < 0$.

We mark transition walls by a line determined by the linear relation, and a decorating arrow attached to the line that indicates the associated sign. For example, in the right hand diagram below, the vertical line corresponds to the relation $a = 2$ with sign $+1$. Note that no $K$-type to the left of that line occurs. The horizontal line corresponds to the relation $b = -2$ with sign $-1$. Presence of $K$-types below this line indicates that the depicted $(\mathfrak{g},K)$-module is reducible.

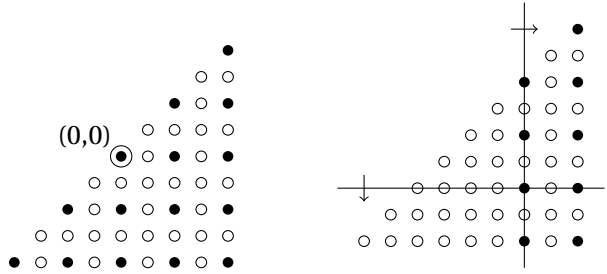

**§2.2 Functions on the upper half space.** We now explain how to pass back and forth between $(\mathfrak{g},K)$-modules and functions $\mathbb{H}^{(2)} \to V(\sigma) \otimes V(\rho)$. Since the group $\Gamma^{(2)}$ and its representation $\rho$ do not enter the construction, we can safely assume that $\rho = \mathbb{1}$. As for the weight, note that any representation $\sigma$ of $\mathrm{GL}_2(\mathbb{C})$ defines a representation of $K$ by restricting along the map $K \to \mathrm{U}_2(\mathbb{R}) \subset \mathrm{GL}_2(\mathbb{C})$, $\begin{pmatrix} a & b \\ -b & a \end{pmatrix} \mapsto a - ib$. When referring to a weight as a $K$-representation we throughout mean this restriction. Given $f : \mathbb{H}^{(2)} \to V(\sigma)$ we consider the function $\mathrm{A}_{\mathbb{R}}(f) : G \to V(\sigma)$ defined by

$$\mathrm{A}_{\mathbb{R}}(f)(g) = \mathrm{A}_{\mathbb{R},\sigma}(f)(g) = \sigma^{-1}\big(j(g, iI_2)\big) f(gi).$$

By contraction of $A_\mathbb{R}(f)$ and $V(\sigma)^\vee$ we obtain a space $\overline{A}_\mathbb{R}(f)$ of functions on $G$, which under right translation yields a $K$-representation that is isomorphic to $\sigma^\vee$. More precisely, for any $v^\vee \in V(\sigma)^\vee$, we have

$$v^\vee\big(A_\mathbb{R}(f)(gk)\big) = v^\vee\big(\sigma^{-1}\big(j(gk,iI_2)\big)f(gi)\big) = v^\vee\big(\sigma^{-1}\big(j(k,iI_2)\big)\sigma^{-1}\big(j(g,iI_2)\big)f(gi)\big)$$
$$= \big(\sigma\big(j(k,iI_2)\big)v^\vee\big)\big(\sigma^{-1}\big(j(g,iI_2)\big)f(gi)\big) = \big(\sigma\big(j(k,iI_2)\big)v^\vee\big)\big(A_\mathbb{R}(f)(g)\big).$$

We refer to the $(\mathfrak{g},K)$-modules generated by the components of $A_\mathbb{R}(f)$ as the $(\mathfrak{g},K)$-module associated with $f$. By the $G$-representation associated with $f$ we mean the representation which is generated by $A_\mathbb{R}(f)$ under right translation. Its Harish-Chandra module is the $(\mathfrak{g},K)$-module that we have associated with $f$.

Vice versa it is possible to extract a function on $\mathbb{H}^{(2)}$ from any vector in a $K$-type of a $(\mathfrak{g},K)$-module that is realized by smooth functions on $G$. For any such $\tilde{f}$ in a representation isomorhic to $\sigma^\vee$, we set

$$A_\mathbb{R}^{-1}(\tilde{f})(\tau) = \sigma\big(j(g_\tau, iI_2)\big)\tilde{f}(g_\tau i).$$

Obviously, we have $A_\mathbb{R}^{-1} A_\mathbb{R}(f) = f$.

We next describe a connection between the differential operators in [Kle+15] and transition of $K$-types in a $(\mathfrak{g},K)$-module associated with $f : \mathbb{H}^{(2)} \to V(\sigma)$. Covariance of the lowering and raising operator and compatibility of $\mathfrak{g}$- and $K$-actions allows us to focus on values of $f$ at $iI_2 \in \mathbb{H}^{(2)}$ and values of $A_\mathbb{R}(f)$ at the identity element $e \in G$. Helgason's formalism [Hel59; Hel62] implies that the action of R at $iI_2$ is given by

$$H^0\big(\mathfrak{k}, \mathfrak{m}^+ \otimes \mathrm{sym}^{\vee 2}\big) \otimes \mathbb{1} \subseteq H^0\big(\mathfrak{k}, \mathfrak{m}^+ \otimes \mathrm{sym}^{\vee 2}\big) \otimes H^0\big(\mathfrak{k}, \sigma \otimes \sigma^\vee\big) \subseteq H^0\big(\mathfrak{k}, \mathfrak{m}^+ \otimes \sigma \otimes (\mathrm{sym}^2 \sigma)^\vee\big),$$

where $\mathbb{1}$ refers to the canonical copy of the trivial representation in $\sigma \otimes \sigma^\vee$. Analogously, the lowering operator L is given by $\mathfrak{k}$-invariants in $\mathfrak{m}^- \otimes (\det^{-2} \mathrm{sym}^2)^\vee$.

Since the action of $\mathfrak{m}$ defines transition of $K$-types, we therefore conclude that the following diagrams commute up to nonzero scalar multiples for all $f \in C^\infty\big(\mathbb{H}^{(2)} \to V(\sigma)\big)$.

$$\begin{array}{ccc}
f & \xrightarrow{A_{\mathbb{R},\sigma}} & \overline{A}_{\mathbb{R},\sigma}(f) \\
R_\sigma \downarrow & & \downarrow \mathfrak{m}^+ \\
R_\sigma f & \xrightarrow{A_{\mathbb{R},\mathrm{sym}^2\sigma}} & \mathfrak{m}^+ \overline{A}_{\mathbb{R},\sigma}(f)
\end{array}
\qquad
\begin{array}{ccc}
f & \xrightarrow{A_{\mathbb{R},\sigma}} & \overline{A}_{\mathbb{R},\sigma}(f) \\
L_\sigma \downarrow & & \downarrow \mathfrak{m}^- \\
L_\sigma f & \xrightarrow{A_{\mathbb{R},\det^{-2}\mathrm{sym}^2\sigma}} & \mathfrak{m}^- \overline{A}_{\mathbb{R},\sigma}(f)
\end{array}$$

The map $A_\mathbb{R}^{-1}$ allows us to deduce vanishing statements with respect to lowering and raising operators from the structure of $(\mathfrak{g},K)$-modules.

**§2.3 Vogan-Zuckerman modules.** We recall from [VZ84] the $(\mathfrak{g},K)$-modules $A_\mathfrak{q}(\lambda)$ that are attached to $\theta$-stable parabolic subalgebras of $\mathfrak{g}$ and central characters $\lambda \in \mathfrak{h}_\mathbb{C}^\vee$ satisfying specific conditions. According to Theorem 5.4 in [VZ84] such representations are uniquely determined by their $K$-types and the central character.

The notion of $A_\mathfrak{q}(\lambda)$ builds up on $\theta$-stable parabolic subalgebras $\mathfrak{q} \subset \mathfrak{g}$, whose definition is revisited on page 56 of [VZ84]. As pointed out there, the name is chosen in an unfortunate way, since not every parabolic subalgebra that is stable under the Cartan involution is neccesary $\theta$-stable in the sense of the following definition. Given $x$ in the Cartan subalgebra of $i\mathfrak{k}_0$, we attach a parabolic subalgebra to $x$ by means of

$$\mathfrak{q} = \mathfrak{q}_x = \mathfrak{l} \oplus \mathfrak{u}, \quad \mathfrak{l} = \ker(\mathrm{ad}(x)), \quad \mathfrak{u} = \text{sum of positive eigenspaces of } \mathrm{ad}(x).$$

We give details for the parabolic subalgebras corresponding to holomorphic and antiholomorpic diskrete series, and to the non-holomorphic Saito-Kurokawa lift. The following choices are suitable elements of $i\mathfrak{k}_0$.

$$x_{\mathrm{hol}} = \mathfrak{h}_\mathfrak{c}, \quad x_{\overline{\mathrm{hol}}} = -\mathfrak{h}_\mathfrak{c}, \quad x_{\mathrm{SK}} = \mathfrak{h}_\mathfrak{k}.$$

For the purpose of studying transitions of $K$-types in $A_{\mathfrak{q}}(\lambda)$ it suffices to determine the intersection of $\mathfrak{u}$ and $\mathfrak{m}$, which is

$$\mathfrak{u}_{\mathrm{hol}} \cap \mathfrak{m} = \mathrm{span}\{\mathfrak{h}_{\mathfrak{m}}^{+}, \mathfrak{e}_{\mathfrak{m}}^{+}, \mathfrak{f}_{\mathfrak{m}}^{+}\}, \quad \mathfrak{u}_{\overline{\mathrm{hol}}} \cap \mathfrak{m} = \mathrm{span}\{\mathfrak{h}_{\mathfrak{m}}^{-}, \mathfrak{e}_{\mathfrak{m}}^{-}, \mathfrak{f}_{\mathfrak{m}}^{-}\}, \quad \mathfrak{u}_{\mathrm{SK}} \cap \mathfrak{m} = \mathrm{span}\{\mathfrak{e}_{\mathfrak{m}}^{+}, \mathfrak{e}_{\mathfrak{m}}^{-}\}. \tag{2.1}$$

In particular, the half-sums of weights in $\mathfrak{u} \cap \mathfrak{m}$ are $\rho_{\mathrm{hol}} = 3\mathfrak{h}_{\mathfrak{c}}$, $\rho_{\overline{\mathrm{hol}}} = -3\mathfrak{h}_{\mathfrak{c}}$, and $\rho_{\mathrm{SK}} = 2\mathfrak{h}_{\mathfrak{k}}$.

**§ 2.4 The Saito-Kurokawa lift.** From Theorem 2.5 of [VZ84], we infer that the minimal $K$-type of $A_{\mathfrak{q}}(\lambda)$ has highest weight $2\rho + \lambda$. In conjunction with Proposition 7.7 of [Miy04], we find that the Harish-Chandra modules associated with $f \in {}^{\mathrm{SK}}\mathrm{S}^{(2)}(\det^{-k/2}\mathrm{sym}^{k} \otimes \rho)$ are isomorphic to $A_{\mathfrak{q}_{\mathrm{SK}}}((k-2)\mathfrak{h}_{\mathfrak{k}}^{\vee}) \otimes V(\rho)$. Note that $(k-2)\mathfrak{h}_{\mathfrak{k}}^{\vee}$ can and will be written as the pair $(k/2-1, 1-k/2)$. When comparing the Harish-Chandra parameter of $A_{\mathfrak{q}}(\lambda)$ to, for example, what Schmidt found on p. 225 of [Sch05], keep in mind that we normalize weights in a way that is convenient from a classical point of view.

**Lemma 2.1.** *Given* $f \in {}^{\mathrm{SK}}\mathrm{S}^{(2)}(\det^{-k/2}\mathrm{sym}^{k})$, *the following vanishing results hold:*

$$\pi_{\mathrm{L},0}\mathrm{L}f = 0, \quad \pi_{\mathrm{L},-}\mathrm{L}f = 0, \quad \pi_{\mathrm{R},0}\mathrm{R}f = 0, \quad \pi_{\mathrm{R},-}\mathrm{R}f = 0.$$

*Proof.* The $K$-type $\overline{\mathrm{A}}_{\mathbb{R}}(f)$ is the minimal one in the Harish-Chandra module associated with $f$. The lowering and raising operators correspond to transitions of $K$-types, as explained in Subsection 2.2. To prove the corollary it therefore suffices to check that the given lowering and raising operators correspond to transitions by $\mathfrak{h}_{\mathfrak{m}}^{-}, \mathfrak{f}_{\mathfrak{m}}^{-}, \mathfrak{h}_{\mathfrak{m}}^{+}$, and $\mathfrak{f}_{\mathfrak{m}}^{+}$. None of them occurs in $\mathfrak{u}_{\mathrm{SK}} \cap \mathfrak{m}$, and therefore the $K$-types corresponding to target weights of the covariant differential operators in question do not occur in $A_{\mathfrak{q}_{\mathrm{SK}}}((k-2)\mathfrak{h}_{k}^{\vee})$. ∎

**§ 2.5 Degenerate generalized principle series.** We determine the composition series of some generalized principal series related to the Eisenstein series that made appearance in [BRR12]. The result that we use is due to Lee [Lee96]. Kudla and Rallis [KR90] previously studied the same principal series, but they did not arrive at the exact conclusion that we need for our purposes. Note also that we could use Muić's study of principal series [Mui09], but then again his approach does not yield exactly what we need.

The next proposition describes induction from the parabolic subgroup

$$P = \left\{ \begin{pmatrix} a & b \\ 0 & d \end{pmatrix} \in G \right\} \quad \text{with} \quad \left| \begin{pmatrix} a & b \\ 0 & d \end{pmatrix} \right| = \det(a),$$

employing unitarily normalization.

**Proposition 2.2 (Kudla-Rallis, Lee).** *Fix a positive, even* $k \in \mathbb{Z}$. *The induced* $(\mathfrak{g}, K)$-*module* $\mathrm{Ind}_{P}^{G}|\cdot|^{k-1/2}$ *has socle series*

$$A_{\mathfrak{q}_{\mathrm{SK}}}(\lambda), \; A_{\mathfrak{q}_{\mathrm{hol}}}(\lambda) \oplus A_{\mathfrak{q}_{\overline{\mathrm{hol}}}}(\lambda), \; V(\lambda),$$

*where* $\lambda = 2(k+1)\mathfrak{h}_{\mathfrak{k}}^{\vee}$ *and* $V(\lambda)$ *is the finite dimensional $G$-representation with infinitesimal character* $\lambda$.

*Proof.* Theorem 5.2 of [Lee96] applies for even $k$. We recall some of the notation of Section 5 in [Lee96]. In Lee's notation we have $n = 2$, $m = 1$, and therefore $r = 1$, and $\alpha = -1 - k$. He defines subquotients $L_{p,q}$ of $\mathrm{Ind}_{P}^{G}|\cdot|^{k-\frac{1}{2}}$ with $1 \leq p, q \leq 2$ by

$$X_{1}^{1} = \{(a,b) \in 2\Lambda : a < k\}, \qquad X_{2}^{1} = \{(a,b) \in 2\Lambda : a \geq k\},$$
$$Y_{1}^{1} = \{(a,b) \in 2\Lambda : b \leq -k\}, \qquad Y_{2}^{1} = \{(a,b) \in 2\Lambda : b > -k\}.$$

The definition of $L_{p,q}$ comprises four cases

$$L_{1,1} = X_{1}^{1} \cap Y_{1}^{1}, \quad L_{1,2} = X_{1}^{1} \cap Y_{2}^{1}, \quad L_{2,1} = X_{2}^{1} \cap Y_{1}^{1}, \text{ and } \quad L_{2,2} = X_{2}^{1} \cap Y_{2}^{1}.$$

Theorem 5.2 of [Lee96] tells us that the socle series of $\mathrm{Ind}_P^G|\cdot|^{k-\frac{1}{2}}$ is $L_{2,1}$, $L_{2,2} \oplus L_{1,1}$, $L_{1,2}$.

All modules are cohomological by inspection of their Harish-Chandra parameters, since $L_{1,2}$ is finite dimensional. To identify the module $L_{2,1}$, we can therefore employ Theorem 5.3 in [VZ84]. Specifically, the infinitesimal character is $2(k+1)\mathfrak{h}_{\mathfrak{c}}^{\vee}$, which under the action of the Weil group is equivalent to $2(k+1)\mathfrak{h}_{\mathfrak{k}}^{\vee}$. This shows that assumption (b) of Theorem 5.3 in [VZ84] is satisfied. The ease verification of assumptions (a) and (c), first notice that $K$-types in $\mathrm{Ind}_P^G|\cdot|^{k-\frac{1}{2}}$ occur with multiplicity at most one by what is explained in Section 2 of [Lee96]. More precisely, we have

$$\left(\mathrm{Ind}_P^G|\cdot|^{k-\frac{1}{2}}\right)_K = \bigoplus_{a,b \in 2\Lambda} V_{K,(a,b)},$$

where $V_{K,(a,b)}$ is the finite dimensional $K$-representation with highest weight $(a,b)$. With this at hand, consider the following diagram $\mathrm{Ind}_P^G|\cdot|^{k-\frac{1}{2}}$.

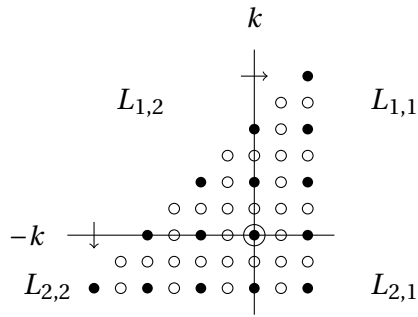

The minimal $K$-type of $L_{2,1}$ is encircled. Inspection of minimal $K$-types and possible transitions of $K$-types verifies assumptions (a) and (c) in each of the infinite dimensional cases. Further, $L_{1,2}$ is finite dimensional, and therefore determined by its infinitesimal character. ∎

## 3 Existence of harmonic weak Siegel Maaß forms

The goal of this section is to prove the first main theorem. It is a consequence of the following one.

**Theorem 3.1.** *Fix $g \in {}^{\mathrm{SK}}\mathrm{S}^{(2)}\bigl(\mathrm{det}^{-k/2}\mathrm{sym}^k \otimes \rho\bigr)$. Then there is a function $f : \mathbb{H}^{(2)} \to V(\mathrm{det}^{2-k/2}\mathrm{sym}^{k-2}) \otimes V(\rho)$ with possible singularities along a suitable divisor D such that:*

*(i) We have $\pi_{\mathrm{L},+}\mathrm{L} f = g$, $\pi_{\mathrm{L},0}\mathrm{L} f = 0$, and $\pi_{\mathrm{L},-}\mathrm{L} f = 0$.*

*(ii) For all $\gamma \in \Gamma^{(2)}$ we have $f|_{\mathrm{det}^{2-k/2}\mathrm{sym}^{k-2},\rho} \gamma = f$.*

*(iii) Singularities of $f$ are meromorphic. That is, for every $\tau \in \mathbb{H}^{(2)}$ there is a neighborhood $U$ of $\tau$ and a meromorphic function $f_U$ on $U$ such that $f - f_U$ extends to $U$.*

A function that satisfies the proberties in the previous theorem will be called a harmonic weak Siegel Maaß form. The space of such forms will be denoted by ${}^{\mathrm{SK}}\mathbb{S}^{(2)}\bigl(\mathrm{det}^{-k/2}\mathrm{sym}^k \otimes \rho\bigr)$.

*Remark 3.2.* Note that harmonic weak Siegel Maaß forms defined above are not the only generalizations of harmonic weak Maaß forms for $\mathrm{SL}_2(\mathbb{R})$ that we expect to exist. Natural candidates arise by searching for preimages under R of holomorphic Siegel modular forms. The methods in this paper, however, do not allow us to construct such preimages, since the correct characterization by means of differential operators would be as follows: Given a holomorphic Siegel modular form $g$ of weight $\mathrm{det}^k$, then a harmonic weak Siegel Maaß forms $f$ in a suitable sense should exist which satisfies

$$\pi_{\mathrm{R},-}\mathrm{R} f = g, \quad \pi_{\mathrm{R},0}\mathrm{R} f = 0, \quad \pi_{\mathrm{L},0}\mathrm{L} f = 0.$$

Such a description by differential operators does not fit into the Dolbeault resolution that we employ in this section.

**§3.1 Toroidal compactifications.** Throughout this section, we let $Y_\Gamma = \Gamma \backslash \mathbb{H}^{(2)}$ for congruence subgroups $\Gamma \subset \mathrm{Sp}_2(\mathbb{Z})$. If $\Gamma$ is neat, then $Y_\Gamma$ is a smooth quasi-projective variety. To simplify notation we will abbreviate $Y_\Gamma$ by $Y$ for one fixed choice of neat $\Gamma$. Only this Subsection 3.1 and 3.2 contain references to $Y_\Gamma$ different from $Y$.

Recall from [Ash+75] that there there is a smooth toroidal compactification $X_\Gamma$ of $Y_\Gamma$ if $\Gamma$ is neat. The corresponding compactification of $Y$ is denoted by $X$. By [San96] the fundamental group of $X_\Gamma$ coincides with the one of $Y_\Gamma$. In the next subsection, we will need to pass between various $X_\Gamma$. Therefore, it is important to note that we may choose compactifications in such a way that for $\Gamma' \subseteq \Gamma$ there is the following commutative diagram of smooth maps between compactifications and coverings. This follows from the dominance statement on in the Main Theorem II of [Ash+75] on page 287.

$$\begin{array}{ccc} Y_{\Gamma'} & \hookrightarrow & X_{\Gamma'} \\ \downarrow & & \downarrow \\ Y_\Gamma & \hookrightarrow & X_\Gamma \end{array}$$

**§3.2 Weights and types as vector bundles.** Every representation $\sigma$ of $\mathrm{GL}_2(\mathbb{C})$ and every representation $\rho$ of $\mathrm{Sp}_2(\mathbb{Z})$ together define a vector bundle over $\mathbb{H}^{(2)}$, denoted by $\mathbb{V}_{\sigma,\rho}$. If the image of $(-I_2, -I_4)$ under $\sigma \otimes \rho$ it the identity, then it descends to the quotient $Y_{\Gamma^{(2)}} = \mathrm{Sp}_2(\mathbb{Z}) \backslash \mathbb{H}^{(2)}$. We denote the resulting vector bundle and its pullbacks to $\Gamma \backslash \mathbb{H}^{(2)}$ for any $\Gamma \subseteq \Gamma^{(2)}$ by the same symbol $\mathbb{V}_{\sigma,\rho}$.

The holomorphic vector bundles $\mathbb{V}_{\sigma,\rho}$ extend uniquely to any toroidal compactification of $Y_\Gamma$. To see this for $\mathbb{V}_{\sigma,\mathbb{1}}$, one employs the Hodge bundle and Schur functors, as is explained, for example, on p. 209 of [Bru+08]. As for the $\mathbb{V}_{\mathbb{1},\rho}$ it suffices to pass to $Y_{\ker\rho}$ and pull back a trivial vector bundle to $X_\Gamma$ along the morphisms in Subsection 3.1. Tensoring $\mathbb{V}_{\sigma,\mathbb{1}}$ and $\mathbb{V}_{\mathbb{1},\rho}$ then yields extensions to $X_\Gamma$ of all $\mathbb{V}_{\sigma,\rho}$.

The dual $\mathbb{V}^\vee_{\sigma,\rho}$ of $\mathbb{V}_{\sigma,\rho}$ is the bundle that corresponds to $\sigma^\vee$ and $\overline{\rho}$. In the case of weights, we can make this more concrete: the dual of $\det^k \mathrm{sym}^l$ is $\det^{-k-l} \mathrm{sym}^l$.

**§3.3 Holomorhpic $n$-forms.** We determine the isomorphism classes of vector bundles associated to the sheafs $\mathbb{E}^{0,1}$ and $\mathbb{E}^{0,2}$ of differential $(0,1)$- and $(0,2)$-forms, and of the sheaf $\Omega^3$ of holomorphic 3-forms. At this occasion we also fix notation $\mathscr{O}_X$ for the structure sheaf of $X$, and $\mathscr{O}_X(D)$ for the sheaf of meromorphic functions with divisor bounded from below by $-D$.

It is quickest to first consider 1-forms, which have local basis $d\tau_1$, $dz$, and $d\tau_2$. We use matrix notation for covariant differentials and obtain

$$d\tau = \begin{pmatrix} d\tau_1 & dz \\ dz & d\tau_2 \end{pmatrix}; \quad d(\tau^{-1}) = -\tau^{-1}(d\tau)\tau^{-1}, \quad \mathfrak{Im}(\tau^{-1})^{-1} d(\overline{\tau}^{-1}) \mathfrak{Im}(\tau^{-1})^{-1} = -\tau y^{-1} (d\overline{\tau}) y^{-1} \tau.$$

This implies that the sheaf of holomorphic 1-forms is isomorphic to sections of $\mathbb{V}_{\mathrm{sym}^2,\mathbb{1}}$. Differential $(0,1)$-forms correspond to the vector bundle $\mathbb{V}_{\det^{-2}\mathrm{sym}^2}$. Taking wedge products we get an isomorphism of the sheaf of differential $(0,2)$-forms with $\mathbb{E}^{0,0} \otimes \mathbb{V}_{\det^{-3}\mathrm{sym}^2}$. Holomophic 3-forms correspond to $\mathbb{V}_{\det^3}$.

As a consequence of the above computation, we find that

$$\mathbb{E}^{0,1} \otimes \mathbb{V}_{\det^{2-k/2}\mathrm{sym}^{k-2},\rho} = \mathbb{E}^{0,0} \otimes \left( \mathbb{V}_{\det^{-k/2}\mathrm{sym}^k,\rho} \oplus \mathbb{V}_{\det^{1-k/2}\mathrm{sym}^{k-1},\rho} \oplus \mathbb{V}_{\det^{2-k/2}\mathrm{sym}^{k-2},\rho} \right).$$

This implies immediately the next lemma.

**Lemma 3.3.** *Given $f \in {}^{\mathrm{SK}}\mathrm{S}^{(2)}(\det^{-k/2}\mathrm{sym}^k)$, there is a unique up to scalar multiples global section of the sheaf $\mathbb{E}^{0,1}_X \otimes \mathbb{V}_{\det^{2-k/2}\mathrm{sym}^{k-2},\rho}$ that corresponds to $f$ when pulled back to a section over $\mathbb{H}^{(2)}$.*

**§3.3.1 Differential operators.** We study the differential $\overline{\partial}$ connecting $\mathbb{E}^{0,0} \otimes \mathbb{V}_{\sigma,\rho}$, $\mathbb{E}^{0,1} \otimes \mathbb{V}_{\sigma,\rho}$, and $\mathbb{E}^{0,2} \otimes \mathbb{V}_{\sigma,\rho}$. Sections of these vector bundles correspond to smooth functions on $\mathbb{H}^{(2)}$ that extend to the boundary of $X$ and which transform like Siegel modular forms of type $\rho$, and weight $\sigma$ in the first case, weight $\det^{-2}\mathrm{sym}^2 \sigma$ in the second one, and weight $\det^{-3}\mathrm{sym}^2 \sigma$ in the third one.

**Lemma 3.4.** *For fixed $\sigma$, there are automorphisms $\phi_1$ of $\det^{-2}\mathrm{sym}^2\sigma$ and $\phi_2$ of $\det^{-3}\mathrm{sym}^2\sigma$ such that the following diagram is commutative.*

$$\begin{array}{ccccc}
H^0\big(\mathbb{E}^{0,0} \otimes \mathbb{V}_{\sigma,\rho}\big) & \xrightarrow{\overline{\partial}} & H^0\big(\mathbb{E}^{0,1} \otimes \mathbb{V}_{\sigma,\rho}\big) & \xrightarrow{\overline{\partial}} & H^0\big(\mathbb{E}^{0,2} \otimes \mathbb{V}_{\sigma,\rho}\big) \\
\downarrow & & \downarrow & & \downarrow \\
C^\infty\big(\mathbb{H}^{(2)} \to V(\sigma \otimes \rho)\big) & \xrightarrow{\phi_1 \mathrm{L}} & C^\infty\big(\mathbb{H}^{(2)} \to V(\det^{-2}\mathrm{sym}^2\sigma \otimes \rho)\big) & \xrightarrow{\phi_2 \pi_{\mathrm{L},0} \mathrm{L}} & C^\infty\big(\mathbb{H}^{(2)} \to V(\det^{-3}\mathrm{sym}^2\sigma \otimes \rho)\big)
\end{array}$$

*Proof.* Observe that $\overline{\partial}$ corresponds for functions on $\mathbb{H}^{(2)}$ to a differential operator of order 1 that annihilates holomorphic functions. It is covariant with respect to $\Gamma^{(2)}$, and by approximation it is therefore covariant with respect to $G$. Further, we check that its image separates points by applying it to germs of polynomials. From the uniqueness statement in Proposition 6.9 of [Kle+15] we therefore deduce that is a nonzero multiple of constituents of the lowering operator L. ∎

**Corollary 3.5.** *The following composition of maps is zero.*

$$^{\mathrm{SK}}\mathrm{S}^{(2)}\big(\det^{-k/2}\mathrm{sym}^k\big) \hookrightarrow H^0\big(X, \mathbb{E}^{0,1} \otimes \mathbb{V}_{\det^{2-k/2}\mathrm{sym}^{k-2},\rho}\big) \xrightarrow{\overline{\partial}} H^0\big(X, \mathbb{E}^{0,2} \otimes \mathbb{V}_{\det^{2-k/2}\mathrm{sym}^{k-2},\rho}\big)$$

*Proof.* This follows when combining the previous lemma with Lemma 2.1. ∎

**§3.4 A vanishing statement.**      Before we can proceed to the proof of Theorem 3.1, we need to establish the following vanishing statement.

**Lemma 3.6.** *Let $D$ be the divisor on $X$ of Igusa's Siegel modular form $\chi_{10}$ (cf. [Igu62]). Fix a weight $\sigma$ and a type $\rho$. Then for sufficiently large $n$, we have*

$$\mathrm{H}^1\big(X, \mathbb{V}_{\sigma,\rho} \otimes \mathcal{O}_X(nD)\big) = 0.$$

*Proof.* Serre Duality relates the left hand side to

$$\mathrm{H}^2\big(X, \Omega^3 \otimes \mathbb{V}_{\sigma^\vee,\overline{\rho}} \otimes \mathcal{O}_X(-nD)\big) = \mathrm{H}^2\big(X, \mathbb{V}_{\det^{3-10n}\sigma^\vee,\overline{\rho}}\big).$$

To show that this vanishes, we employ Corollary 5.6 b) in [EV92] (formulated originally in [GR70]). We have to verify that sections of $\mathbb{V}_{\det^{10n-3}\sigma,\rho}$ yield an embedding of $X$ into some projective space.

Note that by Igusa's work, the vector bundle $\mathbb{V}_{\det^k}$ is very ample for sufficiently large $k$. For fixed $\sigma$ and $\rho$, a result of Serre (cf. Theorem 5.17 of [Har77]) implies that for sufficiently large $k$ global sections of $\mathbb{V}_{\det^k\sigma,\rho}$ generate the corresponding sheaf. This establishes the lemma. ∎

**§3.5 A proof of the existence theorem.**      This subsection will be fully occupied by our proof of Theorem 3.1. We start with a sequence of sheafs similar to the one that appears in the proof of Theorem 3.7 of [BF04]. As opposed to [BF04] we have to establish local lifting as an intermediate step.

Let $D$ be the divisor of Igusa's $\chi_{10}$ as in Lemma 3.6, and let $n$ be large enough so that the Lemma's conclusion holds. We start by considering the following exact sheaf sequence over $X$, which can be obtained by tensoring the Dolbeaut resolution of the structure sheaf $\mathcal{O}_X$ with $\mathbb{V}_{\sigma,\rho} \otimes \mathcal{O}_X(nD)$. Details on the Dolbeaut complex can be found in Section I.3.C of [Dem12].

$$\mathbb{V}_{\sigma,\rho} \otimes \mathscr{O}_X(nD) \hookrightarrow \mathbb{E}_X^{0,0} \otimes \mathbb{V}_{\sigma,\rho} \otimes \mathscr{O}_X(nD) \longrightarrow \mathbb{E}_X^{0,1} \otimes \mathbb{V}_{\sigma,\rho} \otimes \mathscr{O}_X(nD) \longrightarrow \cdots$$

$$\cdots \longrightarrow \mathbb{E}_X^{0,2} \otimes \mathbb{V}_{\sigma,\rho} \otimes \mathscr{O}_X(nD) \longrightarrow \mathbb{E}_X^{0,3} \otimes \mathbb{V}_{\sigma,\rho} \otimes \mathscr{O}_X(nD)$$

Given $f \in {}^{\mathrm{SK}}\mathrm{S}^{(2)}\bigl(\det^{-k/2}\mathrm{sym}^k \otimes \rho\bigr)$ we obtain from it a global section in

$$\mathrm{H}^0\bigl(X, \mathbb{E}_X^{0,1} \otimes \mathbb{V}_{\det^{2-k/2}\mathrm{sym}^{k-2},\rho} \otimes \mathscr{O}_X(nD)\bigr) \subseteq \mathrm{H}^0\bigl(X, \mathbb{E}_X^{0,1} \otimes \mathbb{V}_{\det^{2-k/2}\mathrm{sym}^{k-2},\rho}\bigr).$$

Lemma 3.5 says that this section vanishes under $\overline{\partial}$, mapping it to $(0,2)$-forms. In conjunction with exactness, this implies that it lies, locally, in the image of $\overline{\partial}$ mapping from 0-forms to $(0,1)$-forms. We consider the following image sheaf under $\overline{\partial}$:

$$\overline{\partial}\colon \mathbb{E}_X^{0,0} \otimes \mathbb{V}_{\sigma,\rho} \otimes \mathscr{O}_X(nD) \longrightarrow \mathrm{im}_{\overline{\partial}}\bigl(\mathbb{E}_X^{0,1} \otimes \mathbb{V}_{\sigma,\rho} \otimes \mathscr{O}_X(nD)\bigr) \subseteq \mathbb{E}_X^{0,1} \otimes \mathbb{V}_{\sigma,\rho} \otimes \mathscr{O}_X(nD).$$

We next consider the following short exact sequence.

$$\mathbb{V}_{\sigma,\rho} \otimes \mathscr{O}_X(nD) \hookrightarrow \mathbb{E}_X^{0,0} \otimes \mathbb{V}_{\sigma,\rho} \otimes \mathscr{O}_X(nD) \twoheadrightarrow \mathrm{im}_{\overline{\partial}}\bigl(\mathbb{E}_X^{0,1} \otimes \mathbb{V}_{\sigma,\rho} \otimes \mathscr{O}_X(nD)\bigr)$$

Since $\mathbb{E}_X^{0,0}$ is a fine sheaf, its higher cohomology vanishes. We therefore obtain the long exact sequence

$$\mathrm{H}^0\bigl(X, \mathbb{V}_{\sigma,\rho} \otimes \mathscr{O}_X(nD)\bigr) \hookrightarrow \mathrm{H}^0\bigl(X, \mathbb{E}_X^{0,0} \otimes \mathbb{V}_{\sigma,\rho} \otimes \mathscr{O}_X(nD)\bigr) \longrightarrow \mathrm{H}^0\bigl(X, \mathrm{im}_{\overline{\partial}}\bigl(\mathbb{E}_X^{0,1} \otimes \mathbb{V}_{\sigma,\rho} \otimes \mathscr{O}_X(nD)\bigr)\bigr)$$
$$\downarrow$$
$$\mathrm{H}^1\bigl(X, \mathbb{V}_{\sigma,\rho} \otimes \mathscr{O}_X(nD)\bigr)$$

Lemma 3.6 implies that the obstruction space $\mathrm{H}^1\bigl(X, \mathbb{V}_{\sigma,\rho} \otimes \mathscr{O}_X(nD)\bigr)$ vanishes. In particular, there is a global section of $\mathbb{E}_X^{0,0} \otimes \mathbb{V}_{\sigma,\rho} \otimes \mathscr{O}_X(nD)$ that maps under $\overline{\partial}$ to the section corresponding to $g$. Passing back to functions on $\mathbb{H}^{(2)}$ we finish the proof.

**§3.6 Scalar-valued almost harmonic weak Siegel Maaß forms.** From the conditions on $\mathrm{L}f$ we see that the Harish-Chandra module attached to any harmonic weak Siegel Maaß form is the extension of $A_{\mathfrak{q}_{\mathrm{hol}}}(\lambda)$ by $A_{\mathfrak{q}_{\mathrm{SK}}}(\lambda)$ with $\lambda = (k-2)\mathfrak{h}_{\mathfrak{k}}^{\vee}$, which is a quotient of $\mathrm{Ind}_P^G |\cdot|^{(k-1)/2}$ studied in Section 2.5. This implies that there is a nonzero constant $c_k$ that depends only on $k$ such that

$$\bigl(\pi_{\mathrm{L},+}\mathrm{L}\bigr)^{k/2-1}\bigl(\pi_{\mathrm{R},-}\mathrm{R}\bigr)^{k/2-1} f = c_k f \qquad \text{for all} \quad f \in {}^{\mathrm{SK}}\mathbb{S}^{(2)}\bigl(\det^{2-k/2}\mathrm{sym}^{k-2} \otimes \rho\bigr).$$

From this and the characterization of almost meromorphic Siegel modular forms in terms of the lowering operator, we infer Corollary II.

## 4 Meromorphic and non-holomorphic parts

Given the proof of existence of harmonic weak Siegel Maaß forms, which we achieved in the previous section, we now study Fourier expansions. In particular, we find a canonical lift of Fourier coefficients of non-holomorphic Saito-Kurokawa lifts. This allows us to define Siegel mock modular forms as the remaining terms. The canonical lift of Fourier coefficients stems from Eisenstein series that were investigated in [BRR12].

**§4.1 Fourier coefficients of some Eisenstein series.** Recall from Equation (1.2) that non-holomorphic Saito-Kurokawa lifts have Fourier expansions supported on indefinite coefficients. To phrase (1.2), we employed an ad hoc definition for the analytic part of Fourier coefficients. We are going to replace it by another type of Fourier coefficient that we have better control of. For $k > 2$, we consider the classical Eisenstein series

$$E^{(2)}_{k,\frac{1}{2}} = \sum_{\gamma \in \Gamma^{(2)}_\infty \backslash \Gamma^{(2)}} \det(y)^{\frac{1}{2}}|_{\det^k} \gamma. \tag{4.1}$$

**Proposition 4.1.** *For $k > 2$ the Eisenstein series $E^{(2)}_{k,1/2}$ converges.*

*(1) It has Fourier expansion*

$$\sum_{t \in \mathrm{Mat}_2(\mathbb{Q})} c_E(k;t)\,{}^E a^{(2)}_{k,\frac{1}{2}}(ty)e(tx)$$

*where $c_E(k;t) \in \mathbb{C}$. We have $c_E(k;t) = 0$, if $t$ is negative definite.*

*The analytic part of the Fourier coefficients ${}^E a^{(2)}_{k,1/2}(ty)$ is a function that decays exponentially with respect to the sum of the absolute values of the eigenvalues of $ty$.*

*(2) The Harish-Chandra module associated with $E^{(2)}_{k,1/2}$ is an extension of $A_{\mathfrak{q}_{\mathrm{hol}}}(\lambda)$ by $A_{\mathfrak{q}_{\mathrm{SK}}}(\lambda)$ with Harish-Chandra parameter $\lambda = 2(k+1)\mathfrak{h}^\vee_{\mathfrak{k}}$.*

*Proof.* Observe that $E^{(2)}_{k,1/2}$ can be related to the Eisenstein series studied in [BRR12] as follows:

$$E^{(2)}_{k,\frac{1}{2}} = \det(y)^{\frac{1}{2}} \sum_{\gamma \in \Gamma^{(2)}_\infty \backslash \Gamma^{(2)}} 1|_{k+\frac{1}{2},\frac{1}{2}} \gamma = \det(y)^{\frac{1}{2}} \overline{\left( \sum_{\gamma \in \Gamma^{(2)}_\infty \backslash \Gamma^{(2)}} 1|_{\frac{1}{2},(k+1)-\frac{1}{2}} \gamma \right)},$$

where the group $\Gamma^{(2)}_\infty$ and the slash action on the right hand side are defined by

$$\Gamma^{(2)}_\infty = \left\{ \begin{pmatrix} a & b \\ 0 & d \end{pmatrix} \in \Gamma^{(2)} \right\} \quad \text{and} \quad (f|_{k,k'} \gamma)(\tau) = \det(c\tau + d)^{-k} \det(c\bar{\tau} + d)^{-k'} f(g\tau).$$

The slash action on the right hand side of the above equality is the weight $k+1$ skew slash action that made appearance in [BRR12]. For the purpose of reference, we set

$$^{\mathrm{sk}}E^{(2)}_{k+1} = \sum_{\gamma \in \Gamma^{(2)}_\infty \backslash \Gamma^{(2)}} 1|_{\frac{1}{2},(k+1)-\frac{1}{2}} \gamma = \sum_{t \in \mathrm{Mat}_2(\mathbb{Q})} c_{E,\mathrm{sk}}(k+1;t)\,{}^{\mathrm{sk}}a^{(2)}_{k+1}(ty)e(tx).$$

The Fourier expansion of ${}^{\mathrm{sk}}E^{(2)}_{k+1}$ establishes the one of $E^{(2)}_{k,1/2}$. From Theorem 4 in [BRR12], we see that $c_{E,\mathrm{sk}}(k+1;t) = 0$ if $t$ is positive definite. The relation with $E_{k,1/2}$ implies that $c_E(k;t) = 0$ if $t$ is negative definite. Exponential decay of the analytic part of the Fourier coefficients follows from Shimura's study of confluent hypergeometric functions [Shi82].

To obtain the Harish-Chandra module associated with $E^{(2)}_{k,1/2}$ note that this Eisenstein series yields an intertwining operator from $I_k = \mathrm{Ind}^G_P |\cdot|^{k-1/2}$ into $L^2(\Gamma^{(2)} \backslash G)$. The standard reference for Eisenstein series is [Lan76], but [Lan89] is a survey article which classical readers might find more accessible. The induced representation $I_k$ contains all $K$-types with multiplicity at most one, and therefore $A_{\mathbb{R},\det^k}(E^{(2)}_{k,1/2})$ generates the $K$-type $\det^k$ of the image of $I_k$ in $L^2(\Gamma^{(2)} \backslash G)$. We conclude the proof by applying Proposition 2.2. ∎

*Remark 4.2.* The previous proposition determines the Harish-Chandra module associated with ${}^{\mathrm{sk}}E^{(2)}_{k+1}$ defined in the course of its proof. In particular, we see that Harish-Chandra modules associated with index $t$ Fourier coefficients of ${}^{\mathrm{sk}}E^{(2)}_{k+1}$ for negative definite $t$ are antiholomorphic discrete series. In other words, for negative definite $t$ Fourier coefficients of ${}^{\mathrm{sk}}E^{(2)}_{k+1}$ can be expressed as the product of $e(t\bar\tau)$ and a polynomial of degree $k/2$ in the entries of $y^{-1}$. This implies that the Kohnen Limit process studied in [BRR12] can also be defined for Fourier-Jacobi coefficients of negative index.

**§4.2 An analogue to non-holomorphic Eichler integrals.** We are now ready to define the analytic parts of Fourier coefficients of Saito-Kurokawa lifts. For positive weight $k$ with $4\,|\,k$ we set

$$^{\mathrm{SK}}a^{(2)}_{\det^{2-k/2}\mathrm{sym}^{k-2}}(t;y)e(tx) = (\pi_{\mathrm{L},+}\mathrm{L})^{k/2-1}\Big(^{E}a^{(2)}_{\frac{k}{2},\frac{1}{2}}(ty)e(tx)\Big),$$

$$^{\mathrm{SK}}a^{(2)}_{\det^{-k/2}\mathrm{sym}^{k}}(t;y)e(tx) = (\pi_{\mathrm{L},+}\mathrm{L})^{k/2}\ \Big(^{E}a^{(2)}_{\frac{k}{2},\frac{1}{2}}(ty)e(tx)\Big).$$

As a consequence of Proposition 4.1 and the connection of raising and lowering operators with $(\mathfrak{g},\mathrm{K})$-modules discussed in Section 2.2, we find that

$$\pi_{\mathrm{L},0}\mathrm{L}\big(^{\mathrm{SK}}a^{(2)}_{\det^{2-k/2}\mathrm{sym}^{k-2}}(t;y)e(tx)\big) = 0, \quad \pi_{\mathrm{L},-}\mathrm{L}\big(^{\mathrm{SK}}a^{(2)}_{\det^{2-k/2}\mathrm{sym}^{k-2}}(t;y)e(tx)\big) = 0, \quad \text{and} \tag{4.2}$$

$$\pi_{\mathrm{L},+}\mathrm{L}\big(^{\mathrm{SK}}a^{(2)}_{\det^{2-k/2}\mathrm{sym}^{k-2}}(t;y)e(tx)\big) = {}^{\mathrm{SK}}a^{(2)}_{\det^{-k/2}\mathrm{sym}^{k}}(t;y)e(tx).$$

**Proposition 4.3.** *Consider $f \in {}^{\mathrm{SK}}\mathrm{S}^{(2)}\big(\det^{-k/2}\mathrm{sym}^{k}\otimes\rho\big)$ with Fourier expansion*

$$f = \sum_{t\,\mathrm{indefinite}} c(f;t)\,{}^{\mathrm{SK}}a^{(2)}_{\det^{-k/2}\mathrm{sym}^{k}}(t;y)e(tx).$$

*Then the series*

$$f^{*} = \sum_{t\,\mathrm{indefinite}} c(f;t)\,{}^{\mathrm{SK}}a^{(2)}_{\det^{2-k/2}\mathrm{sym}^{k-2}}(t;y)e(tx)$$

*converges locally absolutely and satisfies*

$$\pi_{\mathrm{L},+}\mathrm{L}f^{*} = f,\quad \pi_{\mathrm{L},0}\mathrm{L}f^{*} = 0,\ \text{and}\quad \pi_{\mathrm{L},-}\mathrm{L}f^{*} = 0.$$

*Proof.* Using exponential decay of ${}^{E}a^{(2)}_{k/2,1/2}(ty)$ stated in Proposition 4.1, we find that

$$\sum_{t\,\mathrm{indefinite}} c(f;t)\,{}^{E}a^{(2)}_{\frac{k}{2},\frac{1}{2}}(ty)e(tx)$$

converges locally absolutely. We can therefore interchange summation and application of lowering operators, to obtain $f^{*}$ as the image under the $(k/2-1)$-th power of $\pi_{\mathrm{L},+}\mathrm{L}$. The behavior of $f^{*}$ under lowering operators follows from (4.2). ∎

**§4.3 Meromorphic parts of harmonic weak Siegel Maaß forms.** The analogue of non-holomorphic Eichler integrals allows us to define Siegel mock modular forms. We call a meromorphic function $f : \mathbb{H}^{(2)} \to V(\det^{2-k/2}\mathrm{sym}^{k-2})\otimes V(\rho)$ a Siegel mock modular form of weight $\det^{2-k/2}\mathrm{sym}^{k-2}$ and type $\rho$, if there is $g \in {}^{\mathrm{SK}}\mathrm{S}^{(2)}\big(\det^{-k/2}\mathrm{sym}^{k}\otimes\rho\big)$ such that

$$\forall\gamma\in\Gamma^{(2)}\,:\,(f+g^{*})|_{\det^{2-k/2}\mathrm{sym}^{k-2},\rho}\,\gamma = f+g^{*}.$$

In this case $g$ is called the shadow of $f$.

Conversely, Proposition 4.3 says that any harmonic weak Siegel Maaß form $f \in {}^{\mathrm{SK}}\mathbb{S}^{(2)}\big(\det^{2-k/2}\mathrm{sym}^{k-2}\otimes\rho\big)$ allows for a decomposition into a meromorphic and a non-holomorphic part. More precisely, we have

$$f = f^{+} + f^{-} \quad \text{with}\quad f^{-} = \big(\pi_{\mathrm{L},+}\mathrm{L}f\big)^{*}\ \text{and}\ f^{+} = f - f^{-}. \tag{4.3}$$

In particular, any harmonic weak Siegel Maaß form is determined up to meromorphic Siegel modular forms by its image under L.

Harmonic Weak Siegel Maaß Forms I *M. Westerholt-Raum*

Chalmers tekniska högskola och Göteborgs Universitet, Institutionen för Matematiska vetenskaper, SE-412 96 Göteborg, Sweden

E-mail: martin@raum-brothers.eu

Homepage: http://raum-brothers.eu/martin